\documentclass[12pt]{amsart}
\usepackage{amscd,amssymb}
\usepackage[graph,frame,poly,arc]{xy}  
\usepackage[plainpages,backref,urlcolor=blue]{hyperref}
\usepackage{mathtools}

\theoremstyle{plain}
\newtheorem{thm}[subsection]{Theorem}
\newtheorem{lem}[subsection]{Lemma}
\newtheorem{prop}[subsection]{Proposition}
\newtheorem{cor}[subsection]{Corollary}

\theoremstyle{definition}
\newtheorem{rk}[subsection]{Remark}
\newtheorem{definition}[subsection]{Definition}
\newtheorem{ex}[subsection]{Example}

\numberwithin{equation}{section}
\newcommand{\OO}{{\mathcal O}}

\newcommand{\F}{{\mathcal F}}

\newcommand{\G}{{\mathcal G}}

\newcommand{\D}{{\mathcal D}}

\newcommand{\al}{{\alpha}}

\newcommand{\de}{{\delta}}

\newcommand{\Z}{\mathbb{Z}}

\newcommand{\C}{\mathbb{C}}
\newcommand{\PP}{\mathbb{P}}

\newcommand{\tM}{\widetilde{M}}

\DeclareMathOperator{\defect}{def}

\DeclareMathOperator{\reg}{reg}
\DeclareMathOperator{\indeg}{indeg}

\DeclarePairedDelimiter\floor{\lfloor}{\rfloor}

\begin{document}

\title [On the Castelnuovo-Mumford regularity of curve arrangements]
{On the Castelnuovo-Mumford regularity of curve arrangements}

\author[Alexandru Dimca]{Alexandru Dimca}
\address{Universit\'e C\^ ote d'Azur, CNRS, LJAD, France and Simion Stoilow Institute of Mathematics,
P.O. Box 1-764, RO-014700 Bucharest, Romania}
\email{dimca@unice.fr}


\subjclass[2010]{Primary 14H50; Secondary  13D02}

\keywords{Castelnuovo-Mumford regularity, plane curve, line arrangement, conic-line arrangement, nodal curve}

\begin{abstract}  The Castelnuovo-Mumford regularity of the Jacobian algebra and of the graded module of derivations associated to a general curve arrangement in the complex projective plane are studied. The key result is an addition-deletion type result, similar to results obtained by
H. Schenck, H. Terao, S. Toh\u aneanu and  M. Yoshinaga, but in which no quasi homogeneity assumption is needed. 
\end{abstract}
 
\maketitle


\section{Introduction} 

Let $S=\C[x,y,z]$ be the polynomial ring in three variables $x,y,z$ with complex coefficients, and let $C:f=0$ be a reduced curve of degree $d\geq 3$ in the complex projective plane $\PP^2$. If $f=f_1\cdot \ldots \cdot f_s$ is the factorization of $f$ into a product of irreducible factors,
we set $C_i:f_i=0$ for $i=1, \ldots ,s$ for the irreducible components of $C$. Then we regard $C$ as the curve arrangement
$$C=C_1 \cup \ldots \cup C_s$$
and denote $d=\deg C$ and $d_i= \deg C_i$.
We denote by $J_f$ the Jacobian ideal of $f$, i.e. the homogeneous ideal in $S$ spanned by the partial derivatives $f_x,f_y,f_z$ of $f$, and  by $M(f)=S/J_f$ the corresponding graded quotient ring, called the Jacobian (or Milnor) algebra of $f$.
Consider the graded $S$-module of Jacobian syzygies of $f$ or, equivalently, the module of derivations killing $f$, namely
$$D_0(f)=\{(a,b,c) \in S^3 \ : \ af_x+bf_y+cf_z=0\}.$$

In this paper we study the Castelnuovo-Mumford regularity of the graded $S$-modules $D_0(f)$ and $M(f)$. Recall that to any graded $S$-module $M$, one can associate a coherent sheaf $\tM$ on $\PP^2$.
We say that $\tM$ is $m$-regular if
$$H^1(\PP^2, \tM(m-1))=H^2(\PP^2, \tM(m-2))=0.$$
The minimal $m$ such that $\tM$ is $m$-regular is called the
Castelnuovo-Mumford regularity of $\tM$ and is denoted by 
$\reg \tM$. Finally we set $\reg M= \reg \tM$, see for instance \cite{Sch0,SchT1} and also \cite[Definition 54]{Sch} for an alternative definition.
Note that for a reduced singular plane curve $C$ of degree $d$, the following inequality holds
\begin{equation}
\label{eq1}
\reg D_0(f) \leq 2d-4
\end{equation}
and the equality holds if $C$ has a unique node as its singular set,
see Remark \ref{rkST} below. On the other hand, for a line arrangement $C:f=0$ the much stronger inequality
\begin{equation}
\label{eq2}
\reg D_0(f) \leq d-2
\end{equation}
holds, and equality takes place if $C$ has only double points, see 
\cite[Corollary 3.5]{Sch0}. The proof of this inequality is based on the following addition-deletion type result. With our notation above, assume that $s>1$ and set 
$$C'=C_1 \cup \ldots \cup C_{s-1}: f'=0.$$
Then, when $C$ is a line arrangement, Schenck shows in \cite{Sch0} that the sheaves 
$$E(f)=\widetilde{ D_0(f)} \text{ and } E(f')(-1)=\widetilde{ D_0(f)}(-1)$$
are related by a short exact sequence of sheaves, from which the conclusion is derived. Similar exact sequences in the case when $C$ is a conic-line arrangement having only quasi homogeneous singularities were considered in \cite{SchT1}, where the authors concentrate on the freeness of such arrangements. The more general situation of a curve arrangement having only quasi homogeneous singularities was considered in \cite{STY}, where an upper bound of $\reg E(f)$ in terms of
$\reg E(f')$ and $d_s=\deg C_s$ when $C_s$ is smooth is given, see \cite[Lemma 3.6]{STY}.
These exact sequences were extended to cover the situation when non 
quasi homogeneous singularities occur, see \cite[Theorem 2.3]{POG}, which can be restated as follows.
First we need some notation. For an isolated hypersurface singularity $(X,0)$ we set
$$\epsilon(X,0)=\mu(X,0)-\tau(X,0),$$
where $\mu(X,0)$ (resp. $\tau(X,0)$) is the Milnor (resp. Tjurina) number 
of the singularity $(X,0)$. We recall that $\epsilon(X,0) \geq 0$ and the equality holds if and only if $(X,0)$ is quasi homogeneous, see \cite{KS0}. For the curves $D_1$, $ D_2$ and $D=D_1 \cup D_2$ and a point $q \in D_1 \cap D_2$, we set
$$\epsilon(D_1,D_2)_q=\epsilon(D_1 \cup D_2,q)-\epsilon(D_1,q)$$
and then define
$$\epsilon(D_1,D_2)=\sum_{q \in D_1 \cap D_2}\epsilon(D_1,D_2)_q.$$
Now we can recall our result in \cite[Theorem 2.3]{POG}, modulo a twist by $-1$.
\begin{thm}
\label{thmB}
With the above notation, assume that $s>1$ and $C_s$ is a smooth curve. 
Then there is an exact  sequence of sheaves on $ \PP^2$ given by
$$ 0 \to E(f')(-d_s) \stackrel{f_s} \longrightarrow  E(f) \to i_{2*}\F \to 0$$
 where  $i_s: C_s  \to \PP^2$ is the inclusion and $\F=\OO_{C_s}(D)$ a line bundle  on ${C_s}$ such that
 $$\deg D=2-2g_s-d_s-r-\epsilon(C',C_s),$$
 where $g_s$ is the genus of the smooth curve $C_s$ and $r$ is the number of points in the reduced scheme of $C' \cap C_s$.
\end{thm}
Using this result,  our generalized version of \cite[Lemma 3.6]{STY} is the following.
\begin{thm}
\label{thm1}
With the above notation, assume that $s>1$ and $C_s$ is a smooth curve of degree $d_s$.
Then
$\reg(D_0(f)) \leq m_0$, where
$$m_0=\max \left ( \reg(D_0(f'))+d_s, 2d_s-3 + \floor*{\frac{r+\epsilon(C',C_s)}{d_s}}\right ).$$

\end{thm}
In fact, our result also corrects a minor error in  \cite[Lemma 3.6]{STY}, see Remark \ref{rkS}. The case when $C$ is a line arrangement was
settled in \cite[Theorem 3.4]{Sch0} and was used to prove the inequality \eqref{eq2}. Theorem \ref{thm1} has the following weaker, but much simpler version.
\begin{cor}
\label{corT1}
With the above notation, assume that $s>1$ and $C_s$ is a smooth curve of degree $d_s$.
Then
$$\reg(D_0(f)) \leq \max \left ( \reg(D_0(f'))+d_s, \deg(C') +2d_s-3\right ).$$

\end{cor}

The following result is the analog of the inequality \eqref{eq2} for the curve arrangements with all the irreducible components smooth. 

\begin{thm}
\label{thm2}
Let $C:f=0$ be a curve arrangement in $\PP^2$ with $d=\deg f$ such that 
 the irreducible components $C_i:f_i=0$ of $C$ are smooth curves, say of degree $d_i$, for all $i=1, \ldots,s$. Then
$$\reg D_0(f) \leq d+ \de-3,$$
where $\de=\max (d_i \ : \ i=1,\ldots ,s)$
and the equality holds if $C$ is a nodal curve.
\end{thm}

\begin{cor}
\label{corT2}
Let $C:f=0$ be a conic-line arrangement with $d=\deg f$. Then
$$\reg D_0(f) \leq d-1$$
and the equality holds if $C$ is a  nodal conic-line arrangement containing at least one smooth conic.
\end{cor}

\begin{rk}
\label{rkI}
The Castelnuovo-Mumford regularity $\reg(D_0(f))$ does not enjoy simple semi-continuity properties, see \cite[Remark 5.3]{DIM}. Hence there seems to be  no simple way to show that the maximal value of 
$\reg(D_0(f))$ in a fixed class of curve arrangements is obtained for the nodal curves in this class, as it is the case in \eqref{eq2} and Theorem \ref{thm2}. On the other hand, a line arrangement $C:f=0$ satisfies the equality in \eqref{eq2} if and only if $C$ is not formal, see \cite[Corollary 7.8]{DST}, and hence $C$ enjoys some geometric properties in this situation. One may ask
whether the conic-line arrangements $C$ for which the equality holds in Corollary \ref{corT2} enjoy also some special properties.
\end{rk}

\section{Some preliminaries } 

We say that $C:f=0$ is an {\it $m$-syzygy curve} if  the module $D_0(f)$ is minimally generated by $m$ homogeneous syzygies, say $r_1,r_2,...,r_m$, of degrees $\al_j=\deg r_j$ ordered such that 
\begin{equation}
\label{eq0}
0 \leq \al_1\leq \al_2 \leq ...\leq \al_m.
\end{equation}

We call these degrees $(\al_1, \ldots, \al_m)$ the {\it exponents} of the curve $C$. 
The smallest degree $\al_1$ is sometimes denoted by ${\rm mdr}(f)$ and is called the minimal degree of a Jacobian relation for $f$. 

The curve $C$ is {\it free} when $m=2$, since then  $D_0(f)$ is a free module of rank 2, see for instance \cite{B+,Dmax,DStRIMS,ST}. In this case
$\al_1+\al_2=d-1$.
Moreover, there are two classes of 3-syzygy curves which are intensely studied, since they are in some sense the closest to free curves.
First, we have the  {\it nearly free curves}, introduced in \cite{DStRIMS} and studied in \cite{AD, B+, Dmax,  MaVa} which are 3-syzygy curves satisfying $\al_3=\al_2$ and $\al_1+\al_2=d$. 
Then, we have the  {\it plus-one generated line arrangements} of level $\al_3$, introduced by Takuro Abe in \cite{Abe},
which are 3-syzygy line arrangements satisfying $\al_1+\al_2=d$. In general, a  3-syzygy curve will be called a {\it plus-one generated curve}  if it satisfies $\al_1+\al_2=d$.

Consider the sheafification 
$$E(f):= \widetilde{D_0(f)} $$
of the graded $S$-module $D_0(f)$, which is a rank two vector bundle on $\PP^2$, see \cite{Se} for details. Moreover, recall that
\begin{equation} \label{equa1} 
E(f)=T\langle C \rangle (-1),
\end{equation}
where $T\langle C \rangle $ is the sheaf of logarithmic vector fields along $C$ as considered for instance in \cite{MaVa,Se}.

\begin{rk} \label{rkD1} 
Note that in \cite[Equation (1.1)]{STY} the vector bundle $T\langle C \rangle$ is denoted by $Der(-\log C)$, and hence we have
\begin{equation} \label{equa2} 
Der(-\log C)=T\langle C \rangle=E(f)(1).
\end{equation}
In particular, this implies
\begin{equation} \label{equa3} 
\reg(Der(-\log C))=\reg(T\langle C \rangle)=\reg(E(f))-1.
\end{equation}
On the other hand, in \cite[Corollary 3.5]{Sch0}, the vector bundle
$E(f)$ is denoted by $\D$, and hence here no twist is involved. Similarly, in
\cite[Formula (1)]{SchT1}, the vector bundle
$E(f)$ is denoted by $\D_0$, and hence again no twist is involved.
\end{rk}

We define the submodule of Koszul-type relations $KR(f)$ to be
the submodule in $D_0(f)$ generated by the following 3 obvious relations of degree $d-1$, namely
$$(f_y,-f_x,0), \  (f_z,0,-f_x)  \text{  and  } (0,f_z,-f_y).$$
Finally, consider the quotient module of essential relations
\begin{equation} \label{eqKR2} 
ER(f)= D_0(f)/KR(f).
\end{equation}
Note that $C:f=0$ is smooth if and only if $ER(f)=0$.
Using this module, we define for a singular curve $C:f=0$ the invariant
\begin{equation} \label{eqMDR} 
{\rm mdr}_e(f) =\min\{ r \in \Z:  ER(f)_r \ne 0 \}.
\end{equation}
We have ${\rm mdr}_e(f) ={\rm mdr}(f) $ when ${\rm mdr}(f) <d-1$.

We introduce the following  invariants associated with the curve $C:f=0$.
\begin{definition}
\label{defCTST} For a homogeneous reduced polynomial $f \in S_d$ one defines
\begin{enumerate}
\item [(i)] the {\it coincidence threshold}\index{coincidence threshold ${\rm ct}(f)$}
$${\rm ct}(f)=\max \{q:\dim M(f)_k=\dim M(g)_k \text{ for all } k \leq q\},$$
with $g$  a homogeneous polynomial in $S$ of the same degree $d$ as $f$ and such that $g=0$ is a smooth curve in $\PP^2$.
\item [(ii)] the {\it stability threshold}
$${\rm st}(f)=\min \{q~~:~~\dim M(f)_k=\tau(C) \text{ for all } k \geq q\}.$$
\end{enumerate}
\end{definition}
In particular, for a smooth curve $C:f=0$  one has ${\rm ct}(f)=\infty$ and ${\rm st}(f)=3(d-2)+1$.
It is clear that  for a singular curve $C:f=0$ one has
\begin{equation}
\label{REL}
{\rm ct}(f) ={\rm mdr}_e(f)+d-2.
\end{equation}

These new invariants ${\rm ct}(f)$ and ${\rm st}(f)$ enter into the following result, see \cite[Corollary 1.7]{Dmax}, where $T=3(d-2)$.
\begin{thm}
\label{ct+st}
Let $C:f=0$ be a degree $d$ reduced curve in $\PP^2$. Then $C$ is a free (resp. nearly free) curve if and only if 
$${\rm ct}(f)+ {\rm st}(f)=T \text{ (resp.  }   \  {\rm ct}(f)+ {\rm st}(f)=T+2). $$
In the remaining cases one has ${\rm ct}(f)+ {\rm st}(f)\geq T+3.$
\end{thm}

To state the following result,  we recall some more notation. 
Let $J=J_f$ be the Jacobian ideal of $f$ and $I=I_f$ be its saturation with respect to the maximal ideal $(x,y,z)$. Then the singular subscheme $\Sigma_f$ of the reduced curve $C:f=0$ is the 0-dimensional scheme defined by the ideal $I$ and we consider the following sequence of defects
\begin{equation} 
\label{eqDEF}
\defect _k\Sigma_f=\tau(C)-\dim \frac{S_k}{I_k}.
\end{equation}
With this notation, one has the following result, see \cite[Theorem 1]{Bull13}, where $g$ is as in  Definition \ref{defCTST} (i).

\begin{thm}
\label{linsys}
Let $C:f=0$ be a degree $d$ reduced curve in $\PP^2$. If $\Sigma_f$ denotes its singular locus subscheme, then
$$\dim M(f)_{T-k}= \dim M(g)_{k} +\defect _k\Sigma_f $$
for $0\leq k \leq 2d-5$. In particular, if $\indeg (I_f) \leq d-2$, then
$${\rm st}(f)=T-\indeg(I_f)+1.$$
\end{thm}
The second claim in Theorem \ref{linsys} follows by taking $k=\indeg (I_f)-1 \leq d-3$ in the first claim and using the obvious equality
$M(g)_j=S_j$ for $j<d-1$.

\begin{lem}
\label{lem1}
Let $C:f=0$ be a reduced plane curve of degree $d$. Then
$$ \reg J_f=\reg(D_0(f))+d-2 \text{  and   } \reg(M(f))=\reg(D_0(f))+d-3.$$
\end{lem}
\proof
The first claim follows from the obvious exact sequence
$$0 \to D_0(f) \to S^3 \to J_f(d-1) \to 0.$$
The second claim follows from the obvious exact sequence
$$0 \to J_f \to S \to M(f) \to 0$$
which implies $\reg(M(f))=\reg(J_f)-1.$
\endproof
The next result gives the relation between these invariants, see \cite[Theorem 3.3]{DIM}.
\begin{thm}
\label{inv}
Let $C:f=0$ be a reduced singular plane curve of degree $d$. Then 
the equality $$ \reg(M(f)) = {\rm st}(f)$$ holds if and only if $C:f=0$ is a free curve. Otherwise, one has $$ \reg(M(f)) = {\rm st}(f)-1.$$
\end{thm}
\begin{rk}
\label{rkST}
It was shown in \cite[Theorem 1.5 and Example 4.3 (i)]{Edin} that
for a reduced singular plane curve of degree $d$ one has
$${\rm st}(f) \leq 3(d-2),$$
and that equality holds when $C$ has a unique node as its singular set.
Such a curve is not free when $d \geq 3$.
It follows that for $d \geq 3$ one has
$$ \reg(M(f)) \leq 3d-7 \text{ and }  \reg(D_0(f)) \leq 2d-4$$
with equalities when $C:f=0$ is a uninodal curve.

\end{rk}
\begin{ex}
\label{exST}
(i) If $C:f=0$ is a free curve of degree $d$ with exponents $(\al_1,\al_2)$ with $\al_1 \leq \al_2$, then one has $1 \leq \al_1 \leq (d-1)/2$ and hence
$${\rm ct}(f)=\al_1+d-2, \  \reg(M(f)) = {\rm st}(f)=2(d-2)-\al_1=d-3+\al_2  \text{ and }
\reg(D_0(f))=\al_2.$$
This follows from relation \eqref{REL}, Theorems \ref{ct+st} and \ref{inv} and Lemma \ref{lem1}.

(ii) If $C:f=0$ is a plus-one generated curve with exponents $(\al_1,\al_2,\al_3)$, then $\al_1+\al_2=d$ and one has
$$ {\rm st}(f)=d-2+\al_3$$
see \cite[Proposition 2.1]{3-syz}. It follows as above that
$${\rm ct}(f)=\al_1+d-2, \  \reg(M(f)) = {\rm st}(f)-1=d-4+\al_3  \text{ and }
\reg(D_0(f))=\al_3-1.$$

\end{ex}
We conclude this section with a local result, needed in the proofs in the next section.

\begin{lem}
\label{lem2}
Consider a reduced plane curve singularity $(D_1,0)$ and a smooth germ $(D_2,0)$ which is not an irreducible component of $(D_1,0)$.
Then
$$(D_1,D_2)_0-\epsilon(D_1,D_2)_0-1 \geq 0,$$
where $(D_1,D_2)_0$ denotes the intersection multiplicity of $(D_1,0)$ and $(D_2,0)$.
\end{lem}
\proof
We have
$$\epsilon(D_1,D_2)_0=\mu(D_1 \cup D_2,0) -\tau(D_1 \cup D_2,0)-(\mu(D_1,0) -\tau(D_1,0)).$$
Using the formula
$$\mu(D_1 \cup D_2,0) =(\mu(D_1,0)+\mu( D_2,0)+ 2(D_1,D_2)_0-1,$$
see \cite[Theorem 6.5.1]{CTC}, the claim in Lemma \ref{lem2}  is equivalent to the much simpler inequality
$$\tau(D_1\cup D_2,0) \geq \tau(D_1,0) + (D_1,D_2)_0.$$
Choose local coordinates at $0 \in \C^2$ such that the smooth germ $(D_2,0)$ is given by $u=0$ and the singularity $(D_1,0)$ by $g=0$.
Let $R=\C\{u,v\}$ be the convergent power series local ring with $\C$ coefficients and variables $u$ and $v$.  Then $g\in R$ is reduced and non divisible by $u$. The singularity $(D_1,0)$ has an associated Tjurina algebra
$$T(g)=R/I_g,$$
where $I_g=(g,g_u,g_v)$, such that $\tau(D_1,0)=\dim T(g)$. Similarly
$T(ug)=R/I_{ug}$ with $I_{ug}=(ug,g+ug_u,ug_v)$ and $\tau(D_1 \cup D_2,0)=\dim T(ug)$. We have the following exact sequence
$$0 \to T(g) \to \frac{R}{(g,ug_u,ug_v)} \to  \frac{R}{(g,u)} \to 0,$$
where the second map is multiplication by $u$ and the third map is the obvious projection. To show that the second map is injective, assume that for $h \in R$ we have $uh \in (g,ug_u,ug_v)$. It follows that there are germs $a,b,c \in R$ such that
$$uh=ag+bug_u+cug_v.$$
This equality implies that $a$ is divisible by $u$, which is not a factor of $g$, and hence $h \in I_g$. This exact sequence implies that
$$\dim \frac{R}{(g,ug_u,ug_v)}=\dim T(g)+\dim \frac{R}{(g,u)}=\tau(D_1,0)+(D_1,D_2)_0.$$
Since $I_{ug} \subset (g,ug_u,ug_v)$, we have
$$\tau(D_1\cup D_2,0) \geq \dim \frac{R}{(g,ug_u,ug_v)}$$
and this completes our proof.

\endproof

\section{The proofs of the  main results} 

\subsection{Proof of Theorem \ref{thm1}} 

If we twist the exact sequence of sheaves in Theorem \ref{thmB} by 
$\OO_{\PP^2}(t)$ and take the associated long cohomology sequence we get
$$H^1(\PP^2,E(f')(t-d_s)) \to H^1(\PP^2,E(f)(t)) \to H^1(C_s,\OO_{C_s}(D_t) )\to $$
$$\to H^2(\PP^2,E(f')(t-d_s)) \to H^2(\PP^2,E(f)(t)) \to 0.$$
Here 
$$\deg D_t=\deg D + td_s=(2+t)d_s-d_s^2-r -\epsilon(C',C_s).$$
The vanishing $H^1(\PP^2,E(f')(t-d_s))=0$ takes place for any
$t-d_s \geq \reg(D_0(f'))-1$, hence for any
\begin{equation}
\label{c1}
t \geq \reg(D_0(f'))+d_s-1.
\end{equation}
Next, we have
$$\dim H^1(C_s,\OO_{C_s}(D_t) )=\dim H^(C_s,\OO_{C_s}(K-D_t) ),$$
where $K$ is the canonical divisor of the curve $C_s$.
It follows that
$$\deg (K-D_t)=\deg(K)-\deg(D_t)=2d_s^2-5d_s+r+\epsilon(C',C_s)-td_s,$$
since $\deg K=d_s(d_s-3)$. It follows that $H^1(C_s,\OO_{C_s}(D_t) )=0$ if $\deg(K-D_t)<0$, in other words if
\begin{equation}
\label{c2}
t > 2d_s-5+\frac{r+\epsilon(C',C_s)}{d_s}.
\end{equation}
This strict inequality is easily  seen to be equivalent to the following
non-strict inequality.
\begin{equation}
\label{c3}
t \geq 2d_s-4+ \floor*{\frac{r+\epsilon(C',C_s)}{d_s}}.
\end{equation}
The inequalities \eqref{c1} and \eqref{c3} imply that the integer $m_0$ defined in Theorem \ref{thm1} satisfies $H^1(\PP^2,E(f)(m_0-1))=0$.
To complete the proof of Theorem \ref{thm1} it remains to show that
$H^2(\PP^2,E(f)(m_0-2))=0$. To get this vanishing, we take $t=m_0-2$ in the above exact sequence and see that 
$$H^2(\PP^2,E(f')(m_0-2-d_s))=0.$$
Indeed, one has 
$$m_0-2-d_s \geq (\reg(D_0(f'))+d_s)-2-d_s=\reg(D_0(f'))-2$$
and we know that $H^2(\PP^2,E(f')(\reg(D_0(f'))-2)=0$. In fact, for any coherent sheaf $\F$ of $\PP^2$, the vanishing $H^2(\PP^2,\F(m))=0$
implies the vanishing $H^2(\PP^2,\F(m+1))=0$ as the obvious exact sequence 
$$0 \to \F(m) \to \F(m+1) \to \G \to 0$$
shows. Here the morphism $\F(m) \to \F(m+1)$ is induced by multiplication by a linear form $\ell \in S_1$, $L$ is the line $\ell=0$ and
$\G$ is a coherent sheaf supported on $L$.
This completes our proof.

\begin{rk}
\label{rkS}
Even in the case when all the singularities in the intersection $C' \cap C_s$ are quasi-homogeneous, our result is slightly different from 
\cite[Lemma 3.6]{STY}. First of all, taking into account the twist explained in Remark \ref{rkD1} and equation \eqref{equa3}, Lemma 3.6 in \cite{STY} can be rested as
$$\reg(D_0(f)) \leq \max \left ( \reg(D_0(f'))+d_s, 2d_s-4 + \floor*{\frac{r}{d_s}}\right ).$$
This difference with our Theorem \ref{thm1} comes from the fact that in
\cite{STY} the strict inequality \eqref{c2} is not replaced by the non-strict inequality \eqref{c3}. When $r/d_s$ is an integer and if
$$\reg(D_0(f'))+d_s \leq 2d_s-4 + \frac{r}{d_s}$$
then the claim in \cite[Lemma 3.6]{STY} is false. Such situations really do occur. Indeed, let $C'$ be a free curve of degree $d'$ and exponents
$(\al_1,\al_2)$ with $\al_2 \leq d'-3$. Let $C_s$ be a smooth curve meeting $C'$ transversally in $d_sd'$ points, which are all nodes for $C$. 
Then Example \ref{exST} (i) implies
$$\reg(D_0(f'))+d_s =\al_2+d_s \leq d'-3+d_s \leq 2d_s-4+d'.$$
Hence such examples exist even in the class of line arrangements.
On the other hand, Theorem 3.4 in \cite{Sch0} which covers the case of 
$C$ a line arrangement is correctly stated.
\end{rk}

\subsection{Proof of Corollary   \ref{corT1}} 
Note that using Lemme \ref{lem2}, we have the following
$$r+\epsilon(C',C_s)=\sum_{p\in C' \cap C_s}(1+\epsilon(C',C_s)_p) \leq \sum_{p\in C' \cap C_s} (C',C_s)_p=\deg(C') d_s.$$
This clearly proves Corollary   \ref{corT1}.
\subsection{Proof of Theorem  \ref{thm2} and Corollary   \ref{corT2}} 
We can assume in this proof that $\de=d_1>1$, since the case of line arrangements is clear by \cite{Sch0}. 
Then we have
$$\reg (D_0(f_1)=2d_1-3,$$
using for instance Theorem \ref{inv}. Hence the first claim holds for $k=1$. Now assume that this claim holds for $k=s-1$ and apply Corollary \ref{corT1}. We get
$$\reg (D_0(f)) \leq \max(d+d_1-3, d_1+ \ldots +d_{s-1}+2d_s-3).$$
This inequality yields the first claim for $k=s$ since $d_1 \geq d_s$.

Now we consider the second claim, when $C$ is in addition a nodal curve. Such a curve $C$ cannot be free, see for instance \cite{DPo}. 
Hence the equality $\reg (D_0(f)) = d+d_1-3$
is equivalent to the equality
\begin{equation}
\label{c13}
{\rm st}(f)=2d-5+d_1,
\end{equation}
in view of Lemma \ref{lem1} and Theorem \ref{inv}. Hence it remains to prove the following.

\begin{lem}
\label{lem3}
Let $C:f=0$ be a nodal curve arrangement in $\PP^2$ with $d=\deg f$ such that 
 the irreducible components $C_i:f_i=0$ of $C$ are smooth curves, say of degree $d_i$, for all $i=1, \ldots,s$. Then,
if
 $\de=\max (d_i \ : \ i=1,\ldots ,s)>1,$
one  has the following equalities
$${\rm st}(f)=2d-5+\de \text{  and  } \indeg I_f=d-\de.$$
\end{lem}
\proof
Assume again that $\de=d_1$.
First  we use Theorem \ref{linsys} and see that  \eqref{c13} is equivalent to
$\indeg I_f=d-d_1$, in other words to the two relations
$$I_{f,d-d_1-1}=0 \text{  and  } I_{f,d-d_1} \ne 0.$$
Since $C$ is a nodal curve, then $I_f$ consists of all the polynomials vanishing at all the nodes of $C$. In particular
$$f_2 f_3 \cdots f_s \in I_{f,d-d_1}$$
and hence $I_{f,d-d_1}\ne 0$.

Finally we prove that $I_{f,d-d_1-1}=0$. Let $h \in I_{f,d-d_1-1}$
and assume first that the curve $H:h=0$ is reduced. For any
$1 \leq k \leq s$, we consider the intersection $H \cap C_k$.
Note that on $C_k$ there are exactly $d_k(d-d_k)$ nodes of the curve $C$. The inequality
$$ d_k(d-d_k) >d_k(d-d_1-1)=\deg C_1 \deg H$$
implies that $C_k$ is an irreducible component of $H$ for all
$k=1, \ldots,s$. This is impossible since 
$$\deg C=d >d-d_1-1=\deg H.$$
If the curve $H$ is not reduced, we apply the above argument to the associated reduced curve $H^{red}$ and get again a contradiction since $\deg H^{red}\leq \deg H$.
This completes the proof of Lemma \ref{lem3} and also of Theorem  \ref{thm2}.  
\endproof

Corollary   \ref{corT2} is an obvious consequence of Theorem  \ref{thm2} for $\de=2$.


\begin{thebibliography}{00}


\bibitem{Abe} T. Abe, Plus-one generated and next to free arrangements of hyperplanes, Int. Math. Res. Not., Vol. 2021, Issue 12 (2021), 9233 -- 9261.


\bibitem{AD}  T. Abe, A. Dimca, On the splitting types of bundles of logarithmic vector fields along plane curves, Internat. J. Math. 29 (2018), no. 8, 1850055, 20 pp.



\bibitem{B+} E. Artal Bartolo, L. Gorrochategui, I. Luengo, A. Melle-Hern\' andez,
On some conjectures about free and nearly free divisors, in: {\it Singularities and Computer Algebra, Festschrift for Gert-Martin Greuel on the Occasion of his 70th Birthday}, pp. 1--19, Springer (2017).









\bibitem{Bull13}  A. Dimca, Syzygies of Jacobian ideals and defects of linear systems,
Bull. Math. Soc. Sci. Math. Roumanie Tome 56(104) No. 2 (2013), 191--203.


\bibitem{Dmax}  A. Dimca, Freeness versus maximal global Tjurina number for plane curves, Math. Proc. Cambridge Phil. Soc.  163 (2017), 161--172.


\bibitem{POG}  A. Dimca, On free and plus-one generated curves arising from free curves by addition-deletion of a line, arXiv: 2310.08972.


\bibitem{DIM} A. Dimca, D. Ibadula, A. M\u acinic, Numerical invariants and moduli spaces for line arrangements, Osaka J. Math. 57 (2020), 847--870.




\bibitem{DPo} A. Dimca, P. Pokora, On conic-line arrangements with nodes, tacnodes, and ordinary triple points, Journal of Algebraic Combinatorics
DOI 10.1007/s10801-022-01116-3.



\bibitem{Edin} A. Dimca, G. Sticlaru, Koszul complexes and pole order filtrations, Proc. Edinburg. Math. Soc. 58(2015), 333--354.






\bibitem{DStRIMS} A. Dimca, G. Sticlaru, Free and nearly free curves vs. rational cuspidal plane curves, Publ. RIMS Kyoto Univ. 54 (2018), 163--179.


\bibitem{3-syz} A. Dimca, G. Sticlaru, Plane curves with three syzygies, minimal Tjurina curves, and nearly cuspidal curves, Geometriae Dedicata 207 (2020), 29--49.

\bibitem{DST} M. DiPasquale, J. Sidman, W. Traves, Geometric aspects of the Jacobian of a hyperplane arrangement,
arXiv:2209.04929 

 



 \bibitem{MaVa} S. Marchesi, J. Vall\` es, Nearly free curves and arrangements: a vector bundle point of view, Math. Proc. Cambridge Philos. Soc. 170 (2021), 51--74.
 
  
  
  
  
 
 \bibitem{KS0} K. Saito, Quasihomogene isolierte Singularit\"aten von Hyperfl\"achen, Invent. Math., 14 (1971), 123--142.




\bibitem{Sch0} H. Schenck, Elementary modifications and line configurations in $\PP^2$. Comment. Math.
Helv. 78 (2003), 447--462. 

\bibitem{Sch} H. Schenck, Hyperplane arrangements: computations and conjectures. In Arrangements of
hyperplanes—Sapporo 2009, volume 62 of Adv. Stud. Pure Math., pages 323--358. Math. Soc.
Japan, Tokyo, 2012.


\bibitem{STY} H. Schenck, H. Terao, M. Yoshinaga,
Logarithmic vector fields for curve configurations in $\PP^2$ with quasihomogeneous singularities,
\textit{Math. Res. Lett.}
\textbf{25}: 1977--1992 (2018).


\bibitem{SchT1} H. Schenck, S. Toh\u aneanu, Freeness of conic-line arrangements in $\PP^2$, Comment. Math. Helv. 84 (2009), 235--258.


\bibitem{Se} E. Sernesi,  The local cohomology of the jacobian ring, {Documenta Mathematica},  19 (2014), 541-565. 






\bibitem{ST} A. Simis, S. O. Toh\u aneanu, Homology of homogeneous divisors, Israel J. Math. 200 (2014), 449-487.








\bibitem{CTC} C. T. C. Wall, \textit{Singular Points of Plane Curves}. Cambridge University
Press, 2004.




\end{thebibliography}
\end{document}